\documentclass[a4paper,11 pt]{article}

\usepackage{latexsym}
\usepackage[english]{babel}
\usepackage[utf8]{inputenc}  
\usepackage{amsmath} 
\usepackage{amssymb} 
\usepackage{mathrsfs}
\usepackage{graphicx}
\usepackage{ esint }
\usepackage[usenames]{color}
\usepackage{enumerate}

\def\epsilon{\varepsilon}

\newtheorem{lemm}{Lemma}[section]

\newtheorem{teor}[lemm]{Theorem}

\newtheorem{deff}[lemm]{Definition}
\newtheorem{oss}[lemm]{Remark}

\def\dim{{\bf Proof.}}

\def\a1{$a_{1}$}
\def\r1{$r_{1}$}

\def\1{\{}
\def\2{\}}

\newcommand{\cvd}{\begin{flushright}$\Box$\end{flushright}}

\newcommand{\eq}{\begin{equation}}
\newcommand{\feq}{\end{equation}}
\newcommand{\be}{\begin{equation}}
\newcommand{\ee}{\end{equation}}

\begin{document}
\begin{center}

\vspace{1.cm}

{\Huge
\textbf{Energy-exchange stochastic models for non-equilibrium}}

\vspace{0.5cm}

\begin{center}

{\Large
Chiara Franceschini$^{1}$, Cristian Giardin\`a$^{1}$} 
\vspace{0.2cm}

\textit{$^{1}$ Department of Mathematics,\\
University of Modena and Reggio Emilia}\\
\end{center}

\end{center}

\vspace{0.5cm}

\noindent
Non-equilibrium steady states are subject to intense investigations 
but still poorly understood. For instance, the derivation of Fourier law 
in Hamiltonian systems is a problem that still poses several obstacles. 
In order to investigate non-equilibrium systems,
stochastic models of energy-exchange have been introduced 
and they have been used to identify universal  properties of non-equilibrium.
In these notes, after a brief review of the problem of anomalous transport
in 1-dimensional Hamiltonian systems, some boundary-driven interacting 
random systems are considered and the ``duality approach'' to their rigorous
mathematical treatment is reviewed. 
Duality theory, of which a brief introduction is given, is a powerful technique 
to deal with Markov processes and interacting particle systems. 
The content of these notes is mainly based on the papers \cite{GK,GKR,GKRV}.


\vspace{1.cm}

\noindent
These notes are based on two lectures given by C. Giardin\`a during the training week in the GGI workshop 'Advances in Nonequilibrium Statistical Mechanics: large deviations and long-range correlations, extreme value statistics, anomalous transport and long-range interactions' in Florence, Italy (May 2014). These notes were prepared by Chiara Franceschini.

\newpage

\section{Fourier's law and anomalous transport}
The models we would like to discuss are motivated by 
two fundamental open problems in mathematical statistical physics: 
\begin{enumerate}
\item
deriving phenomenological laws of
non-equilibrium statistical physics for a microscopic
interacting model;
\item
understanding the general structure
and properties of the  probability measure
describing  a system in a non-equilibrium steady state.
\end{enumerate}

\noindent
In this opening section we briefly discuss the first problem,
while we defer the second to the study of the specific models
discussed below. We focus on the Fourier's law
\begin{equation}
\label{fourier2}
\langle J \rangle =\kappa\cdot \nabla T
\end{equation}
describing the heat flow across a metal bar when we heat the
bar from one side and cool it from the other. The law
bears its name to  J.B.J. Fourier \cite{fourier} who 
discovered it in 1822. In  \eqref{fourier2} 
 the left hand side $\langle J \rangle$  is the 
average energy current per unit time and per unit of surface, whereas in the right side
$\kappa$ is the material's conductivity and $\nabla T$ is the 
temperature gradient. One would like to derive 
this phenomenological law (and the linear heat equation describing
the diffusive heat spreading) from a simplified, yet realistic, mathematical model.
The simplest setting is obtained by considering
a one-dimensional system with $L$ sites
coupled to two external reservoirs imposing 
temperatures $T_{l}$ and $T_{r}$ at the extremes.
A crucial distinction emerges as those two parameters
are varied:
\begin{itemize}
\item[-]
If $T_{l}=T_{r}$ then in the long time limit the system reaches
an equilibrium state described by the Boltzmann-Gibbs probability 
distribution.
\item[-]
If $T_{l}\neq T_{r}$ then a non-equilibrium state arises.
We will be interested in those cases in which
a stationary measure sets in the long time limit. 
Such invariant measure will be called in the sequel
the non-equilibrium  probability measure.
\end{itemize}

\subsection{Hamiltonian models}
\noindent
In the basic setting described above one considers a bulk part 
and a boundary contribution. If one starts from 
the assumption that the micro-world evolution is described 
by Newton equations, then the bulk part of the model  
consists of particles 
whose dynamics is encoded in the Hamiltonian
\begin{equation}
\label{hami}
H_L(q,p) = \sum_{i=1}^{L} \left(\dfrac{p_{i}^{2}}{2} + U(q_{i})\right) + \sum_{i=1}^{L-1}\; V(q_{i+1}-q_{i})
\end{equation}
Here the particles are assumed to have unit mass, position and momentum 
coordinates are denoted by
$(q,p)=(q_1,\ldots,q_L,p_1,\ldots,p_L)$. The Hamiltonian has a local contribution 
(including kinetic energy and a site potential $U$) and a nearest-neighbor interaction 
potential $V$.
For the boundaries, one option is to model the reservoirs by adding 
to the velocities of the first and the last particles an Ornstein-Uhlenbeck 
process that fixes the temperatures via a fluctuation-dissipation mechanism\footnote{It is also possible to work with deterministic thermostats 
(e.g. Nos\'e-Hoover or isokinetic), however this setting will not be discussed 
here (see \cite{reviewLLP} for more on this)}.
Therefore the full equations of motion read
\begin{equation}
\left\{
\begin{array}{l}
dq_{i}  = p_{i}dt\\
\\
dp_{i}  =  -\dfrac{\partial H_L}{\partial q_{i}}dt 
+ \delta_{i,1}(-\gamma p_{1}+\sqrt{2 \gamma T_{l}} \ dW^{(l)})
+ \delta_{i,L}(-\gamma p_{L}+\sqrt{2 \gamma T_{r}} \ dW^{(r)})
\end{array}
\right.
\nonumber
\end{equation}
where $W^{(l)}$ and $W^{(r)}$ are two independent standard Brownian motion
and $0<\gamma<\infty$ is a parameter tuning the coupling to reservoirs.
The microscopic definition of the observables appearing in the 
Fourier law \eqref{fourier2} follows from the discretization
of the continuity equation 
$$
\frac{\partial \rho_i}{\partial t} = - (j_{i} - j_{i-1})
$$
where $\rho_i$ denotes the energy density  at site $i$
and $j_{i}$ is the  energy current across the 
bound $(i,i+1)$.
As discussed in \cite{reviewLLP}, this leads to the definition of 
the current in the bulk, i.e. for $i\in\{1,\ldots,L-1\}$
\begin{equation}
\label{loc-current}
j_{i}:= - \dfrac{1}{2} \left( p_i + p_{i+1}\right)\dfrac{\partial H_L}{\partial q_{i}} 
\end{equation}
and the average  current for a system of size $L$ is given by
\begin{equation}
\label{current}
\langle J^{(L)} \rangle : = \dfrac{1}{L-1} \sum_{i=1}^{L-1} \langle j_{i} \rangle
\end{equation}
where $\langle \cdot \rangle$ denotes expectation with respect to 
the stationary non-equilibrium probability measure. 
The details of the computation leading to the definition in \eqref{loc-current}
and \eqref{current}
can be found in \cite{reviewLLP}.
As for the temperature the standard definition is given
by twice the average kinetic energy
 \footnote{Other definitions of the temperature are possible, cfr \cite{R, GL}.},
yielding
\begin{equation}
\label{temp}
T_{i}:=\langle p_{i}^{2}\rangle
\end{equation}
Combining together \eqref{current} and \eqref{temp} and assuming validity
of Fourier's law \eqref{fourier2} with linear stationary energy profiles, 
one obtains a definition of the conductivity
for a system of size $L$
$$
\kappa_{L}: = \langle J^{(L)} \rangle \dfrac{ L}{T_{l} -T_{r}}
$$
and considering the thermodynamic limit one has  the
definition of the system conductivity 
$$
\kappa
:=\lim_{L\rightarrow \infty} \kappa_{L}\;.
$$
\subsection{Stylized properties of Fourier's law in $1d$ Hamiltonian models}
With reference to the generic Hamiltonian \eqref{hami}, the following
picture emerges from numerical and analytical studies of several models.
\begin{enumerate}[i)]
\item For harmonic oscillators with $U(x) = 0$ and 
$V(x)=\dfrac{1}{2}x^{2}$ Lebowitz, Lieb and Rieder \cite{LLR} 
proved that $\kappa_{L}\simeq L$.
The result is rooted in the fact that the $L$ degrees of freedom
can be decoupled into normal modes, which transport ballistically 
the heat from one side to the other. Furthermore  \cite{LLR} proves also 
that the non-equilibrium invariant measure is a multivariate Gaussian 
measure. 

\item For non-linear oscillator chains with a non-vanishing on site 
potential (i.e. $U(x)\neq 0$) a finite conductivity 
is found \cite{reviewLLP}. The reason is that the on-site potential acts as
a source of scattering among the normal modes. However this
case is believed to be quite unrealistic.

\item For non-linear oscillator chains with translation invariant interaction
(i.e. $U(x)=0$) one typically finds $\kappa_{L}\sim L^{\alpha}$ 
with $0< \alpha < 1$. This is the phenomenon of {\em anomalous transport},
whose origin has been linked to the presence of additional conserved
quantities (momentum, besides energy) for the bulk dynamics.
The result is supported both by numerical analysis of several models
(a much studied case is the FPU-$\beta$ model with a quartic potential 
$V(x)=\dfrac{1}{2}x^{2}+\dfrac{1}{4}\beta x^{4}$) as well as by analytical 
studies based on mode-coupling theory \cite{reviewLLP} and nonlinear 
fluctuating hydrodynamics \cite{vB,S}. From the numerical experiments the value
of the exponent $\alpha \in [0.3,0.5]$.

\item An exceptional case is given by the rotor model with potential 
$U(x)=0$, $V(x)=1-cos(x)$. In this case, despite conservation of
momentum, the finite volume conductivity scales ad 
$k_{L}=a+\frac{b}{L}$, yielding a finite conductivity in the
thermodynamic limit. This was found by two (independent) numerical 
studies of the works \cite{GLPV, GS}.

\item 
Recently it has been claimed \cite{ZZWZ} that, despite
the divergence observed in numerical simulations of finite
system sizes, non-linear oscillator chains with $U(x) =0$ 
and asymmetric potential $V(x)\neq -V(x)$ 
have finite asymptotic thermal conductivity 
at low temperatures and anomalous transport at
high temperatures.
The claim has been contradicted in \cite{DDN}.
\end{enumerate}

\subsection{Stochastic models}

The use of stochastic models to model the bulk system is a
further simplifying assumption. 
In this approach the Hamiltonian dynamics is replaced by a stochastic
evolution and exact solutions can be obtained.
The first model of this type was introduced by
Kipnis, Marchioro and Presutti \cite{KMP}  in 1982.
They considered a model in which the energy 
is uniformly redistributed among nearest
neighbor particles, thus providing an efficient mechanism
of energy transport across the extended system.
They proved the validity of Fourier law and introduced 
the duality approach that is the core of these lectures. 
More recent works include \cite{BBO}, 
where the case of harmonic oscillators
with an energy conserving stochastic noise
has been studied. Remarkably they
find that if the noise is only energy-conserving 
then the conductivity is finite. On the contrary,
if the noise conserves both energy and momentum
then $\kappa_{L}\sim \sqrt{L}$, thus strengthening 
the claim that conservation of momentum leads
(in general) to anomalous transport.

\subsection{From Hamiltonian to stochastic}
\begin{figure}[h!]
\centering
\includegraphics[width=7.cm,angle=-90]
{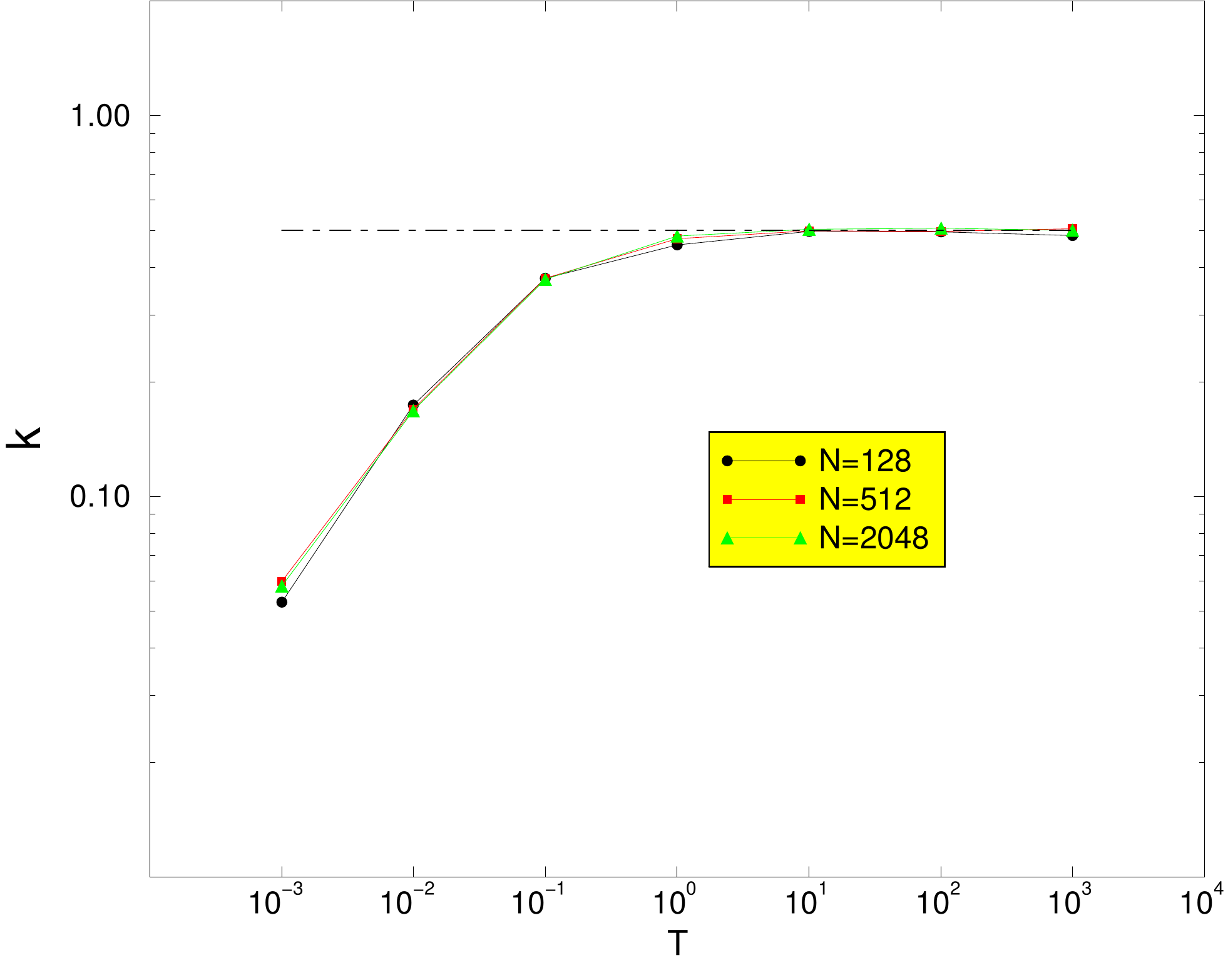} 
\caption{\small{Heat conductivity as a function 
of the temperature for the model \eqref{pippo} 
with reservoirs and for a particular choice of the vector potential $A$. 
Picture taken from \cite{GK}.}} 
\label{fig-kappa}
\end{figure}

We conclude this introduction by considering a model that
has been introduced in \cite{GK}, it serves as
a minimalistic model to go from Hamiltonian to
stochastic dynamics, in the sense that in the 
high-energy limit the deterministic dynamics is 
well-approximated by a stochastic dynamics.
Consider the Hamiltonian
\begin{equation}
\label{pippo}
H_L(q,p)= \sum_{i=1}^{L} \dfrac{(p_{i}-A_{i}(q))^{2}}{2}
\end{equation}
where $A(q)=\bigl( A_{1}(q),\ldots, A_{L}(q) \bigl)$ 
is a generalized ''vector'' potential in $\mathbb{R}^{L}$.
In this model there is a non-trivial distinction between momentum 
and velocities. 
The equations of motion read
\begin{equation*}
\begin{cases}
\dfrac{dq_{i}}{dt}=v_{i} \\ \\
\dfrac{dv_{i}}{dt}= \sum_{j}B_{ij}v_{j}
\end{cases}
\end{equation*}
where the ``magnetic fields" are obtained from the generalized
potential as
$$
B_{i,j}(q)=\dfrac{\partial A_{i}(q)}{\partial q_{j}}-\dfrac{\partial A_{j}(q)}{\partial q_{i}}
$$ 
Since the fields form and antisymmetric matrix, i.e. $B_{i,j} = -B_{j,i}$,  
the total energy is conserved: 
\begin{equation*}
\dfrac{d}{dt} \sum_{i} \dfrac{v_{i}^{2}}{2}=\sum_{i} v_{i} \dfrac{dv_{i}}{dt}= \sum_{i,j} B_{ij}v_{i}v_{j}=0.
\end{equation*}
The model \eqref{pippo}, coupled to reservoirs, has been studied using numerical simulations in \cite{GK} for several choices of the vector potential and the thermal conductivity measured for different system sizes and different temperature. 
From those studies it has been found that the system has always a finite 
conductivity. Moreover, as the temperature is increased, the conductivity approaches
an asymptotic constant value (see figure \ref{fig-kappa}).\\

The findings of the numerical study of the model suggested
that in the high-energy limit the deterministic dynamics could be 
substituted with a random evolution.
The stochastic model that will be discussed in these lectures 
is indeed obtained by replacing the deterministic magnetic 
fields $B_{i,j}$ with a family of random fields given by 
independent Brownian motions \cite{GK}.

\section{The Brownian Momentum Process (BMP)}\label{tbmp}

\begin{deff}[BMP, stochastic differential equations]
For a graph $G_L=(V_L,E_L)$ where $V_L$ is the vertex set of cardinality $L\in \mathbb{N}$ and $E_L$ is the edges set, the Brownian Momentum Process $\left( X(t) \right) _{t\geq 0}=\left( X_{i}(t) \right) _{i \in V_L ,\,  t\geq 0}$ is a diffusion process that takes value in $\mathbb{R}^{L}$ and satisfies the SDE (in Stratonovich sense) 
$$
dX_{i}= \sum_{j: \; (i,j)\in E_L}X_{j} \circ dW_{ij}(t)
$$
where for $1\le i<j \le L$ the $W_{ij}(t)$ are a family of independent standard Brownian motion and for $1\ge i>j \ge L$ we set $W_{ij}(t)=-W_{ji}(t)$. Moreover $W_{ii}(t)=0$ for
$i\in V_L$.
\end{deff}
\begin{oss} 
The $X(t)$ vector collects the velocities of $L$ particles at time $t>0$. 
For any graph $G_L$ it follows from the definition above
that the total kinetic energy is conserved
$$
\dfrac{d}{dt}\left( \frac{1}{2} \sum_{i\in V_L} X_{i}^{2}(t) \right)=0 
$$
\end{oss}
An alternative definitions of a continuous time Markov process is obtained
by specifying its generator.
\begin{deff}[BMP, generator]
The generator of the BMP process on the graph $G_L = (V_l,E_L)$
is the second order differential operator $\mathscr{L}^{BMP}$
defined on smooth functions
$f: \mathbb{R}^L \to \mathbb{R}$ as
\begin{equation}
\label{L-BEP}
\mathscr{L}^{BMP}f(x)=\sum_{(i,j)\in E_L} \mathscr{L}^{BMP}_{ij}f(x)=\sum_{(i,j)\in E} \left( x_{i}\frac{\partial}{\partial x_{j}}-x_{j}\frac{\partial}{\partial x_{i}} \right) ^{2} f(x)
\end{equation}
\end{deff}
\begin{oss} 
We remind that the generator $\mathscr{L}$ of a Markov process gives the infinitesimal evolution of the expectations
of observables. Namely, for a measurable function $f$ 
\begin{equation*}
\mathscr{L}f(x)= \lim_{h \to 0} \dfrac{\mathbb{E}_{x}\left[ f(X_{h})\right] -f(x)}{h}=\dfrac{d}{dt}\mathbb{E}_x\left[ f(X_{t})\right] \bigg|_{t=0}  
\end{equation*}
where $\mathbb{E}_{x}\left[ f(X_{t})\right] = \mathbb{E}\left[ f(X_{t}) \; | \; X_0 = x \right]$.
As a consequence  
$$
\dfrac{d}{dt} \mathbb{E}\left[ f(X_{t})\right] =\mathbb{E}\left[ \mathscr{L}f(X_{t})\right]\;.
$$
By exponentiating the generator one (formally) obtains the semigroup $S_t = \exp[t \mathscr{L}]$. Namely, assuming the Hille-Yoshida theorem \cite{Ligg} can be applied, 
one has
$$
\mathbb{E}_{x}\left[ f(X_{t})\right] = S_t f(x).
$$
\end{oss}
The last characterization of a Markov stochastic process is obtained
by the forward Kolmogorov equation (also called Fokker-Planck equation
in the context of diffusions). For the forward Kolmogorov equation to exists
it is required that the process is absolutely 
continuous with respect to the Lebesgue measure.
\begin{deff}[BMP, Fokker-Planck equation]
The probability densi\-ty function $p(x,t) : \mathbb{R}^L \times \mathbb{R}_+ \to \mathbb{R}$ of the BMP process at time $t>0$ on the graph $G_L = (V_L,E_L)$ started from the density
$g$ at time $t=0$
is given by
\begin{equation}
\begin{cases}
\dfrac{d}{dt}p(x,t)=\bigl( (\mathscr{L}^{BEP})^{*} p \bigl)(x,t)  \\
p(x,0)=g(x)
\end{cases}
\end{equation}
where $(\mathscr{L}^{BEP})^{*}$ is the adjoint in $L^2(dx)$ of $\mathscr{L}^{BEP}$ in \eqref{L-BEP}. It turns out that $(\mathscr{L}^{BEP})^{*} = \mathscr{L}^{BEP}$.
\end{deff}
\begin{oss}
From the Fokker-Planck equation it is easy to check that pro\-duct measure with 
marginal given by ${\cal N}(0,T)$, i.e. centered Gaussian with variance $T$, are 
invariant. It is believed that the family of such product measure, labeled by the
variance $T>0$, exhaust all the ergodic measures \cite{GRV}. In other words it should
be possible to write any other invariant measure as a convex combination
of those product measures.
\end{oss}

\subsection{Brownian Momentum Process with reservoirs}
In order to make connection with the heat transport discussed in the
first section we specify to a graph $G_L$, which is given by the
one-dimensional lattice of $L$ sites with edges between nearest
neighbor vertices. Moreover we add to the first and last site
two reservoirs at temperatures  $T_{l}$ and $T_{r}$ modeled
as Ornstein-Uhlenbeck processes. 
\begin{deff}[BMP with reservoirs]
The generator of the BMP process with reservoirs is 
\begin{equation}
\mathscr{L}^{BMP, \ res}= \sum_{i=1}^{L-1} \left( x_{i}\dfrac{\partial}{\partial x_{j}}-x_{j}\dfrac{\partial}{\partial x_{i}} \right) ^{2}-x_{1}\dfrac{\partial}{\partial x_{1}} + T_{l}\dfrac{\partial ^{2}}{\partial x_{1}^{2}} -x_{L}\dfrac{\partial}{\partial x_{L}} + T_{r}\dfrac{\partial ^{2}}{\partial x_{L}^{2}}
\end{equation}
\end{deff}
If $T_{l}=T_{r}=T$ then a unique stationary measure is selected, namely the product of centered 
Gaussians with variance $T$ with probability density function
\begin{equation*}
p(x)=\prod_{i=1}^{L}\dfrac{1}{\sqrt{2\pi T}}e^{-\frac{x_{i}^{2}}{2T}}.
\end{equation*}
The question we would like to address is
what can be said about the situation in which $T_{l}\neq T_{r}$.
We will see that a characterization of the non-equilibrium
invariant measure can be achieved by using stochastic duality.

\section{Duality theory}
Duality theory is a powerful technique to deal with stochastic processes. In a nutshell, the main idea behind duality is to study a given process making use of a simpler dual one. In particular the connection between the two processes occurs on a set of so-called duality functions.

\subsection{Duality}

\begin{deff}[Duality]
Let $X=(X(t))_{t\geq 0}$ and $Y=(Y(t))_{t\geq 0}$ be two Markov processes with state spaces $\Omega$, respectively $\Omega^{dual}$, and generator $ \mathscr{L}$, respectively $\mathscr{L}^{dual}$.
We say that $X$ is dual to $Y$ with duality function $D:\Omega\times \Omega^{dual}\longmapsto \mathbb{R} $ if 
\begin{equation}\label{D}
\mathbb{E}_{x}[D(X(t),y)]=\mathbb{E}_{y}[D(x,Y(t))] \ ,
\end{equation}
for all $(x,y) \in \Omega \times \Omega^{dual}$ and $t\geq 0$.
\end{deff}

\begin{oss}
Markov processes $X$ and $Y$ are dual on function $D$ iff
\begin{equation}\label{dua}
[\mathscr{L}D(\cdot,y)](x)=[\mathscr{L}^{dual}D(x,\cdot)](y)
\end{equation}
Indeed, informally, we have the following series of equalities
\begin{equation*}
\mathbb{E}_{x}[D(X_{t},y)]=\bigl( e^{t\mathscr{L}}D(\cdot,y) \bigl)(x) = \bigl( e^{t\mathscr{L}^{dual}}D(x,\cdot) \bigl)(y) =\mathbb{E}_{y}[D(x,Y_{t})]
\end{equation*}
\end{oss}

\subsection{Algebraic approach to duality theory}

How to find a dual process? How to find duality-function? 
The following scheme has been put forward in a series of works  \cite{GKR,GKRV}
\begin{itemize}
\item Duality arises as a change of representation of an abstract operator
that belongs to a Lie-algebra. 
\item Duality functions correspond to the 
intertwiners between the two representations.
\end{itemize}
We shall illustrate this approach by considering the case of
Brownian Momentum Process and its underlying $SU(1,1)$ 
Lie algebra structure. 
The existence of a dual process often allows to simplify the
analysis of the process at hand. In the context of boundary
driven non-equilibrium systems, further simplifications
take place. A list of consequences of duality theory includes:
\begin{enumerate}
\item {\em From continuous to discrete:} interacting diffusions can be studied via interacting particles systems.
\item {\em From reservoirs to absorbing boundaries:} stationary state of boun\-dary-dri\-ven 
processes with reservoirs can be fully characterized by dual processes with absorbing boundaries.
\item {\em From many to few:} $n$-point correlation functions of a system of site $L$ can be studied using $n$ dual walkers.
\end{enumerate}

\section{The Symmetric Inclusion Process (SIP)}

\begin{deff}[SIP(m), generator]
Given the graph $G_L=(E_L,V_L)$ the Symmetric Inclusion Process with parameter $m\in\mathbb{R}_+$ is the continuous time Markov chain $(\eta(t))_{t>0} = (\eta_i(t))_{i\in V_L,t\ge 0}$ with state space $\mathbb{N}^{L}$ with 
generator
\begin{align*}
& \mathscr{L}^{SIP(m)} f(\eta) = \sum_{(i,j) \in E_L}  \mathscr{L}_{i,j}^{SIP(m)} f(\eta) = \\ & \sum_{(i,j) \in E_L} \eta_{i} \left( \dfrac{m}{2} + \eta_{j} \right) \left[ f(\eta^{i,j})-f(\eta) \right]+\eta_{j} \left( \dfrac{m}{2} + \eta_{i} \right) \left[ f(\eta^{j,i})-f(\eta) \right]
\end{align*}
where the vector $\eta^{i,j}$ is obtained form the configuration $\eta$ by moving one particle
from site $i$ to site $j$, i.e.
\begin{equation*}
\eta^{i,j}=(\eta_{1},\ldots ,\eta_{i}-1,\ldots, \eta_{j}+1, \ldots,  \eta_{L})
\end{equation*}
\end{deff}
\begin{oss}
The SIP process is an interacting particles system with attractive interactions.
The stationary reversible measures are product with marginals 
Neg Bin$ \left( \dfrac{m}{2},p\right)$, i.e.
\begin{equation*}
\mathbb{P}(\eta_{1}=n_{1},\ldots,\eta_{V}=n_{V})= \prod_{i \in V} p^{n_{i}}(1-p)^{\frac{m}{2}} \binom{n_{i}+\frac{m}{2}-1}{\frac{m}{2}-1}
\end{equation*}
We remind that Neg Bin$ \left( \dfrac{m}{2},p\right)$ represents the number of trials before the $(\dfrac{m}{2})^{th}$ failure in a sequence of independent Bernoulli trials with success probability $p$.
\end{oss}

\subsection{SIP with absorbing boundaries}
\begin{deff}[SIP(1) with absorbing boundaries]
We add two extra sites, i.e. the configurations are $\eta=(\eta_{0},\eta_{1},\ldots , \eta_{L}, \eta_{L+1})$. The generator of Symmetric Inclusion Process with absorbing boundaries and parameter $m=1$ is 
\begin{align*}
\mathscr{L}^{SIP(1),abs} f(\eta) = \sum_{i=1}^{L-1} \;\; 
& \eta_{i} \left( \dfrac{1}{2} + \eta_{i+1} \right) \left[ f(\eta^{i,i+1})-f(\eta) \right]+ \\ 
&  \eta_{i+1} \left( \dfrac{1}{2} + \eta_{i} \right) \left[ f(\eta^{i+1,i})-f(\eta) \right]+ \\ & 
 \dfrac{ \eta_{1}}{2} \left[ f( \eta^{1,0})-f( \eta) \right] + \dfrac{ \eta_{L}}{2} \left[ f(\eta^{L,L+1})-f(\eta) \right]
\end{align*}
\end{deff}

\subsection{Duality between BMP with reservoirs and  SIP(1) with absorbing boundaries}

The main result, which allows studying the non-equilibrium invariant measure
of the Brownian Momentum Process, is contained in the following theorem.
\begin{teor}[Duality BMP/SIP(1) ]
The BMP with reservoirs is dual to the  SIP(1) with absorbing boundaries on duality function
\begin{equation*}
D(x, \eta) = T_{l}^{\eta_{0}} \left(  \prod_{i=1}^{L} \dfrac{x_{i}^{2\eta_{i}}}{(2\eta_{i}-1)!!}  \right)  T_{r}^{\eta_{L+1}} 
\end{equation*}
\end{teor}
\dim \ \  From equation \eqref{dua} the proof is a consequence of the identity 
\begin{equation*}
(\mathscr{L}^{BMP,res}D(\cdot,\eta))(x)=(\mathscr{L}^{SIP(1), abs}D(x,\cdot))(\eta)
\end{equation*}
The above equation can be verified by a direct explicit computation.
\cvd

\section{Intermezzo:  SU(1,1) algebra}

We recall that if $A$ and $B$ are operators working on a common domain, the commutator of $A$ and $B$ is $[A,B]:=AB-BA$.
\begin{deff}[SU(1,1) algebra]
For a graph $G_L = (V_L,E_L)$, we define the SU(1,1) algebra as the algebra of elements generated by $\lbrace K_{i}^{+}, K_{i}^{-}, K_{i}^{0} \rbrace_{i \in V_L}$ satisfying the commutation relations 
\begin{equation}\label{SU}
[K^{-}_{i},K^{+}_{j}]=2K_{i}^{0}\delta_{ij}, \ \ \ \ \ \ \ \ \ \ [K^{0}_{i},K^{\pm}_{j}]=\pm K_{i}^{\pm}\delta_{ij}
\end{equation}
\end{deff}

\noindent
The duality between BMP and SIP process can be seen as a change of representation of the abstract operator $\textbf{L}$ that is a linear combination of the generators of the SU(1,1) algebra
\begin{align*}
\textbf{L} & =\sum_{i=1}^{L-1} \left(  K_{i}^{+}K_{i+1}^{-}+K_{i}^{-}K_{i+1}^{+}-2K_{i}^{0}K_{i+1}^{0}+\dfrac{1}{8} \right)+
\\ & - \left( 2K_{1}^{0}+\dfrac{1}{2} + T_{l}2K_{1}^{-}  \right) - \left( 2K_{L}^{0}+\dfrac{1}{2} + T_{r}2K_{1}^{-}  \right)
\end{align*}
There exists a representation in terms of differential operators where the
generators of the SU(1,1) algebra are given by
\begin{equation*}
K^{+}_{i}=\frac{1}{2} x_{i}^{2}, \ \ \ \ \ \ \ K^{-}_{i}=\frac{1}{2}\partial_{i}^{2}, \ \ \ \ \ \ \ K^{0}_{i}=\frac{1}{4}(x_{i}\partial_{i}+\partial_{i}x_{i}).
\end{equation*}
In this case one finds that $\textbf{L}=\mathscr{L}^{BMP,res}$.
Another discrete representation is obtained in terms of infinite dimensional matrices by writing
\begin{equation*}
\mathcal{K}_{i}^{+}|\eta_i\rangle=(\eta_{i}+\frac{1}{2})| \eta_i+1\rangle \ \ \ \ \   \mathcal{K}_{i}^{-}|\eta_i\rangle=\eta_{i}| \eta_i-1\rangle  \ \ \ \ \ \mathcal{K}_{i}^{0}|\eta_i\rangle=(\eta_{i}+\frac{1}{4})| \eta_i\rangle
\end{equation*}
where $| \eta_i \rangle$ denotes the vector with all components equal to $0$ except the $i^{th}$
component equal to $1$. 
In this case one finds that $\textbf{L}=\mathscr{L}^{SIP(1),abs}$.
Duality functions are the intertwiner between the two representations
and they are obtained by imposing relation (\ref{dua}) for all the algebra 
generators. In particular,
\begin{equation*}
K_{i}^{+}d(\cdot,\eta_{i})(x_{i})=\mathcal{K}_{i}^{+}d(x_{i},\cdot)(\eta_{i})
\end{equation*} 
leads straight to 
\begin{equation*}
d(x_{i},\eta_{i})=\dfrac{x_{i}^{2 \eta_{i}}}{(2 \eta_{i}-1)!!}.
\end{equation*}

\section{Correlation functions in the stationary state}
Next theorem shows that the moments of energy (i.e. square of the velocity) of the BMP process with
reservoirs can be characterize via duality.
\begin{teor}[Moments of BMP]
Let $\eta$ be a configuration of the SIP(1) process with absorbing boundaries and denote
by  $|\eta|=\sum_{i=1}^{L} \eta_{i}$ the total number of SIP dual walkers. Let 
\be
p_{\eta}(a,b)= \mathbb{P}(\eta_{0}(\infty)=a, \ \eta_{L+1}(\infty)=b  \; | \;\eta(0)=\eta )\;.
\ee 
be the probability that $a$ SIP walkers are absorbed at site $0$ and $b$ of them
are absorbed at site $L+1$. Then denoting by $\langle \cdot \rangle_{L}$ the expectation in the non- equilibrium stationary state of the BMP process with reservoirs  one has 
\begin{equation}
\langle D(X,\eta)\rangle_{L}= \sum_{a,b \ : \ a+b=|\eta|} T^{a}_{l}T^{b}_{r}p_{\eta}(a,b)
\end{equation}
\end{teor}
\dim \\
Consider the  BMP process started from the initial measure  $\nu$. Then
\begin{align*}
\langle D(X,\eta)\rangle_{L} & 
=
\lim_{t\rightarrow \infty} \int \mathbb{E}_{x}[D(X(t),\eta)] \nu(dx) \\ & 
= \lim_{t\rightarrow \infty} \int \mathbb{E}_{\eta}[D(x,\eta(t))] \nu(dx)  \\ &  
=\int \lim_{t\rightarrow \infty} \mathbb{E}_{\eta}[D(x,\eta(t))] \nu(dx) 
\\ & = \int \nu(dx) \mathbb{E}_{\eta}(T^{\eta_{0}(\infty)}_{l}T^{\eta_{L+1}(\infty)}_{r}) \\ & 
=\sum_{a,b \ : \ a+b= \eta} p_{\eta}(a,b)T_{l}^{a}T_{r}^{b}
\end{align*}
\cvd
\begin{flushleft}
\textbf{Example 1: temperature profile}
\end{flushleft}
In particular it follows that if $\eta=(0,\ldots,0,1,0,\ldots,0)$ is a null vector with $1$ in the $i^{th}$ position, then $D(X,\eta)=X^{2}_{i}$. Denoting by $(R(t))_{t>0}$ a continuous time
symmetric random walker jumping at rate 1/2 and absorbed at the boundaries $\{0,L+1\}$
also with rate 1/2 one has
\begin{align*}
\langle X^{2}_{i}\rangle 
&= T_{l} \;p_{\eta}(1,0)+T_{r} \;p_{\eta}(0,1) \\ 
& = T_{l} \;\mathbb{P}(R(\infty)=0|R(0)=i)+T_{r} \; \mathbb{P}(R(\infty)=L+1|R(0)=i)  \\ 
& =T_{l}\left( 1- \dfrac{i}{L+1}\right)+T_{r}\left( \dfrac{i}{L+1}\right)  \\ & =T_{l}+\dfrac{T_{r}-T_{l}}{L+1}i
\end{align*}
\begin{flushleft}
\textbf{Example 2: energy covariance}
\end{flushleft}
If $\eta=(0,\ldots,0,1,0,\ldots,0,1,0, \ldots,0)$ is a null vector with $1$ in the $i^{th}$ position and in the $j^{th}$ position, then $D(X,\eta)=X^{2}_{i}X^{2}_{j}$. Therefore
\begin{align*}
\langle X^{2}_{i}X^{2}_{j}\rangle 
&= T_{l}^2 \;p_{\eta}(2,0) + T_{r}^2 \;p_{\eta}(0,2)  + T_l T_r  \; (p_{\eta}(1,0) +  p_{\eta}(0,1)) \\
& = T_{l}^{2} \mathbb{P}_{ij}(R_{1}(\infty)=R_{2}(\infty)=0)+
T_{r}^{2}\mathbb{P}_{ij}(R_{1}(\infty)=R_{2}(\infty)=L+1) \\
& +
T_{l}T_{r} [\mathbb{P}_{ij}(R_{1}(\infty)=0, R_{2}(\infty)=L+1) + 
\mathbb{P}_{ij}(R_{1}(\infty)=L+1, R_{2}(\infty)=0)] 
\end{align*}
where $(R_1(t),R_2(t))_{t>0}$ are two continuous time
SIP walkers absorbed at the boundaries $\{0,L+1\}$.
A computation gives \cite{GKR}
$$
\langle X^{2}_{i}X^{2}_{j}\rangle  - \langle X^{2}_{i}\rangle \langle X^{2}_{j}\rangle 
=
\frac{2i (L+1-j)}{(L+3)(L+1)^2} (T_l - T_r)^2
$$

\section{Redistribution model}
The example of duality between BMP and SIP(1) process can be generalized in
several ways. For instance one can define energy redistribution jump process
by considering instantaneous thermalization limit of the BMP process. 

\subsection{The Brownian Energy Process BEP(m)}
It is convenient to start from a ladder-graph with $m$ copies of the one dimensional
lattice of $L$ sites. The Brownian Momentum Process on such graph has generator
\begin{equation}
\mathscr{L}^{BMP(m)}=\sum_{i=1}^{L-1} \sum_{\alpha,\beta=1}^{m} \left( x_{i,\alpha}\partial_{i+1, \beta}- x_{i+1,\beta} \partial_{i, \alpha}  \right) ^{2}
\end{equation}
Let $Z_{i}(t)= \sum_{\alpha=1}^{m} X_{i,\alpha}^{2}(t)$ be the energy at site $i$ at time $t\ge 0$. Then one can check that $(Z(t))_{t\ge 0} = (Z_i(t)_{i\in V_L, t\ge 0})$ defines a Markov process called the Brownian Energy Process with parameter $m\in \mathbb{R}$ (BEP($m$)) having generator
\begin{equation}
\label{bep-m}
\mathscr{L}^{BEP(m)}=\sum_{i=1}^{L-1} z_{i}z_{i+1} (\partial_{i}-\partial_{i+1})^{2}-\dfrac{m}{2}(z_{i}-z_{i+1})(\partial_{i}-\partial_{i+1})
\end{equation}
For the BEP($m$) process the total energy $\sum_{i=1}^{L} Z_{i}(t)$ is conserved.
It can be proved that the stationary reversible measure are product with marginals $Gamma \left(\dfrac{m}{2},\theta \right)$, i.e. with probability density 
\begin{equation}
p(z)= \prod_{i=1}^{L} \dfrac{z_{i}^{\frac{m}{2}-1}e^{- \frac{z}{\theta}}}
{\Gamma\left( \frac{m}{2} \right)  \theta^{\frac{m}{2}}}
\end{equation}
Exploiting a change of representation of the SU(1,1) algebra from
\begin{equation}\label{cont}
K_{i}^{+}=z_{i} \ \ \ \ \ \ \ \ \  K_{i}^{-}=z_{i}\partial_{i}^{2}+\dfrac{m}{2}\partial_{i} \ \ \ \ \ \ \ \ \  K_{i}^{0}=z_{i} \partial_{i} + \dfrac{m}{4}
\end{equation}
to
\begin{equation}\label{discr}
K_{i}^{+}| \eta_{i}\rangle = \biggr( \eta_{i}+\dfrac{m}{2} \biggr) | \eta_{i}+1\rangle  \ \ \ \ \ \ \  K_{i}^{-}|\eta_{i}\rangle=\eta_{i}|\eta_{i}-1\rangle \ \ \ \ \ \ \  K_{i}^{0}|\eta_{i}=\biggr( \eta_{i}+ \dfrac{m}{4} \biggr) | \eta_{i}\rangle
\end{equation}
one deduces that the BEP(m) process admits a dual given by  the SIP(m) process with generator\begin{eqnarray}
\label{sip-m}
\mathscr{L}^{SIP(m)} f(\eta) & =  & \sum_{i} \eta_{i} \left( \dfrac{m}{2} + \eta_{i+1} \right) \left[ f(\eta^{i,i+1})-f(\eta) \right]  \nonumber\\
&& \quad + \;\eta_{i+1} \left( \dfrac{m}{2} + \eta_{i} \right) \left[ f(\eta^{i+1,i})-f(\eta) \right]
\end{eqnarray}

\noindent
The following result is a generalization of the duality relation between BMP and SIP models.
\begin{teor}[Duality BEP(m)/SIP(m)]
The process with generator $\mathscr{L}^{BEP(m)}$ in \eqref{bep-m} and the process with generator $ \mathscr{L}^{SIP(m)} $ in \eqref{sip-m} are dual with duality function 
\begin{equation*}
D(z,\eta)=\prod_{i=1}^{L} z_{i}^{\eta_{i}} \dfrac{\Gamma \left( \frac{m}{2}\right) }{2^{\eta_{i}}  \Gamma \left( \frac{m}{2}+ \eta_{i} \right)}.
\end{equation*}
\end{teor}

\subsection{Instantaneous thermalization limit}

The idea behind instantaneous thermalization limit is to imagine that
the process evolves through jumps and at each jump it immediately
thermalize with respect to the invariant measure. For the BEP(m) 
this leads to the following definition for the generator of the
instantaneous thermalization limit on the edge $(i,i+1)$
\begin{align*}
&\mathscr{L}_{i,i+1}^{IT, BEP(m)}f\bigl(z_{i},z_{i+1}\bigl):=\lim_{t\rightarrow \infty} \left(\left( e^{tL^{BEP(m)}}-1 \right) f\right)\bigl(z_{i},z_{i+1}\bigl)=\\ & \int \rho^{m}(z_{i}',z'_{i+1}|z_{i}'+z'_{i+1}= z_{i}+z_{i+1})\biggr[ f\bigl(z_{i}',z'_{i+1}\bigl)-f\bigl(z_{i},z_{i+1}\bigl)\biggr] dz_{i}'dz_{i+1}'
\end{align*}
where $\rho^{(m)}$ denotes the probability density of the stationary state of the BEP(m) process on two sites. Since we know that for two independent random variables $(Z_1,Z_2)$ with distribution $Gamma \left(\dfrac{m}{2},\theta \right)$ the ratio $P = Z_1/(Z_1+Z_2)$ is distributed like $Beta \left(\frac{m}{2},\frac{m}{2}\right)$ then we can also write
\begin{align*}
&\mathscr{L}_{i,i+1}^{IT, BEP(m)}f\bigl(z_{i},z_{i+1}\bigl)
=\\ &
\int_{0}^{1}dp \ \nu^{(m)}(p) \biggr[f\bigl(p(z_{i}+z_{i+1}),(1-p)(z_{i}+z_{i+1})\bigl)-f\bigl(z_{i}+z_{i+1}\bigl) \biggr]
\end{align*}
where $\nu^{(m)}$ denotes the probability density of $Beta \left(\frac{m}{2},\frac{m}{2}\right)$
distribution. In particular, since $Beta (1,1)$ coincides with the uniform distribution, then
for $m=2$ the generator of the KMP process \cite{KMP} is recovered, i.e.
\be
\mathscr{L}_{i,i+1}^{IT, BEP(2)}f\bigl(z_{i},z_{i+1}\bigl) =
\mathscr{L}_{i,i+1}^{KMP}f\bigl(z_{i},z_{i+1}\bigl) 
\ee
where
\be
\nonumber
\mathscr{L}_{i,i+1}^{KMP}f\bigl(z_{i},z_{i+1}\bigl) =
\int_{0}^{1}dp \biggr[f\bigl(p(z_{i}+z_{i+1}),(1-p)(z_{i}+z_{i+1})\bigl)-f\bigl(z_{i}+z_{i+1}\bigl)\biggr]
\ee

\section{Duality and multiple conservation laws}

In this last chapter we present an example of a diffusion process that conserves (in the bulk) its total energy and momentum
(see also \cite{BBO,BO}). We investigate its duality relations and, last, we infer a general theorem about duality and change of coordinates. The results discussed in this section are taken from \cite{CF}.

\subsection{A diffusion process with conservation of energy and momentum}

The basic process we examine is a Markov process defined by its generator as follows.
\begin{deff}[of the process]
Consider a diffusion process $ (\textbf{X}_{t})_{t\geq 0}=\left( X_{t},Y_{t},Z_{t} \right)_{t\geq 0} $ taking values in $\mathbb{R}^{3}$. The vector $(x,y,z)$ represents the momentum associated with three unit mass particles $\lbrace 1,2,3 \rbrace $ freely moving in a physical volume $V$.
Thus, up to irrelevant constant, the total momentum of the system is $P=x+y+z$ and the total (kinetic) energy is $E=x^{2}+y^{2}+z^{2}$.
The generator of the process is
\begin{equation}\label{ll}
\mathscr{L} = \bigr[(x\partial_{y}-y\partial_{x})+(y\partial_{z}-z\partial_{y})+(z\partial_{x}-x\partial_{z}) \bigr]^2
\end{equation}
where we shorthand $\partial_{i}=\frac{\partial}{\partial i}$ with $i \in \lbrace  x,y,z \rbrace$.
$\mathscr{L}$ acts on twice differentiable functions $f:\mathbb{R}^{3}\rightarrow \mathbb{R}$.
\end{deff}
A distinguishing property of this process regards the conservation of both total energy $E$ and total momentum $P$. This property can be easily proved by letting act the generator on those functions. It is easy to see that they are both zero.
The same result can also be achieved via the stochastic differential equations (in Itô sense) of $ \mathscr{L} $ associated to generator \eqref{ll}, which are
\begin{align} \label{ito} \begin{cases}
&dx_{t}=(-2x_{t}+y_{t}+z_{t})dt+(y_{t}-z_{t})dB_{t} \\ 
&dy_{t}=(x_{t}-2y_{t}+z_{t})dt+(z_{t}-x_{t})dB_{t}  \\
&dz_{t}=(x_{t}+y_{t}-2z_{t})dt+(x_{t}-y_{t})dB_{t}.  
\end{cases}
\end{align}
By inspection
it is possible to find that the total momentum $P$ is conserved
\begin{align*}
d(x_{t}+y_{t}+z_{t})=&(-2x_{t}+y_{t}+z_{t})dt+(y_{t}-z_{t})dB_{t}+(x_{t}-2y_{t}+z_{t})dt+ \\ & (z_{t}-x_{t})dB_{t}+ (x_{t}+y_{t}-2z_{t})dt+(x_{t}-y_{t})dB_{t}=0
\end{align*}
and, by making use of the Itô's formula, the total energy of the process $E$ satisfies 
\begin{align*}
d(x_{t}^{2}+y_{t}^{2}+z_{t}^{2}) & = \bigl[ 2x_{t}(-2x_{t}+y_{t}+z_{t})+2y_{t}(x_{t}-2y_{t}+z_{t})+2z_{t}(x_{t}+y_{t}-2z_{t}) \\ &
 + 2(y_{t}-z_{t})^{2} + 2(z_{t}-x_{t})^{2}+ 2(x_{t}-y_{t})^{2} \bigl]dt  \\ & + 2\bigl[ x_{t}y_{t}-x_{t}z_{t}+y_{t}z_{t}-y_{t}x_{t}+z_{t}x_{t}-z_{t}y_{t} \bigl]dB_{t} \\ &
 =0.
\end{align*}
Before going any further, we want to highlight the geometric aspects of our problem.
Generator \eqref{ll} can be viewed as the result of three different rotations.
To be more specific, $(x\partial_{y}-y\partial_{x}) $ generates the rotation across the $z$-axis, $(y\partial_{z}-z\partial_{y}) $ generates the rotation across the $x$-axis and $(z\partial_{x}-x\partial_{z}) $ generates the rotation across the $y$-axis, so it turns out that generator \eqref{ll} represents the rotation around the $(1,1,1)$-axis, which is orthogonal to the plane of equation $x+y+z=constant$.

\noindent
As a consequence of the conservation of both total energy and total momentum, the motion takes place in the $1-$dimensional manifold (i.e. a circle) given by the intersection of the sphere $x^{2}+y^{2}+z^{2}=E$ and the plane $x+y+z=P$ orthogonal to the rotation axis just mentioned.

\begin{oss}[Extended system]
It is easy to define a one-dimensional system with $L$ sites coupled with two external reservoirs at different temperatures $T_{l}$ and $T_{r}$ that conserve in the bulk both energy and momentum.
The generator of this process is given by the sum of each generator of the type of \eqref{ll} for the bulk part, where the contribution of the two reservoirs has been added:
\begin{equation}
\mathscr{L} = \mathscr{L}_{l}+ \sum_{i=1}^{L-2} \mathscr{L}_{i,i+1,i+2} +\mathscr{L}_{r}.
\end{equation}
\end{oss}

\subsection{Duality results}

One might wonder about the existence of a dual process in the setting of multiple conservations laws.
In order to find a duality relation we make a change of coordinates that simplifies the expression of the generator \eqref{ll}.
This is achieved by a rotation that maps the $(1,1,1)$ axis to the $z'$ axis of the new coordinates system.
In other words, we need to find matrix $R$ such that 
$$
\left (
\begin{array}{ccc}
\frac{1}{\sqrt{3}} \\
\frac{1}{\sqrt{3}} \\
\frac{1}{\sqrt{3}} \\
\end{array}
\right ) = R \left (
\begin{array}{ccc}
0 \\
0 \\
1 \\
\end{array}
\right ).$$
Matrix $R$ describes the rotation that moves a frame $Ox'y'z'$, initially aligned with $Oxyz$, into a new orientation in which the $Oz'$ axis is brought into the $(1,1,1)$ axis.
This matrix is easily found through Euler angles
\begin{equation}\label{matrixR}
R = \left (
\begin{array}{cc c}
-\frac{\sqrt{2}}{2} cos\varphi -\frac{\sqrt{2}}{2\sqrt{3}} sin\varphi & \frac{\sqrt{2}}{2} sin\varphi -\frac{\sqrt{2}}{2\sqrt{3}} cos\varphi & \frac{1}{\sqrt{3}} \\
\frac{\sqrt{2}}{2} cos\varphi -\frac{\sqrt{2}}{2\sqrt{3}} sin\varphi & -\frac{\sqrt{2}}{2} sin\varphi -\frac{\sqrt{2}}{2\sqrt{3}} cos\varphi & \frac{1}{\sqrt{3}} \\
\frac{\sqrt{2}}{\sqrt{3}}sin\varphi & \frac{\sqrt{2}}{\sqrt{3}}cos\varphi & \frac{1}{\sqrt{3}} \\
\end{array}
\right ).
\end{equation}
This change of coordinates let us find the generator of the process that conserves total momentum and energy (and we call it $\mathscr{L}_{3}$ to highlight the fact that $3$ sites are involved) as function of $x'$, $y'$ and $z'$
\begin{equation}\label{elle3}
\mathscr{L}_{3} = \bigl[(x \partial _{y}-y \partial _{x})+(y \partial _{z}-z \partial _{y})+(z \partial_{z}-x \partial _{z}) \bigl]^2 = 3 \  [x'\partial_{y'} - y'\partial_{x'}]^2.
\end{equation}
To make it clearer we call the latter equation above $\mathscr{L'}_{3}$, i.e. the generator of the Brownian Momentum process with $2$ sites $(x \partial_{y} - y \partial_{x})^{2}$ for which duality is verified, as shown in \cite{GKR}.\\
Thereby the duality function is obtained by the knowledge of 
\begin{equation*}
D'(x',y';n_{1},n_{2})=d(x';n_{1})d(y';n_{2})=\dfrac{x'^{2n_{1}}y'^{2n_{2}}}{(2n_{1}-1)!!(2n_{2}-1)!!}
\end{equation*}
which is the duality function of the process of generator $\mathscr{L'}_{3}$. \\ 
Thanks to matrix $R$ in \eqref{matrixR}, it is possible to write $x'$ and $y'$ as functions of $x$, $y$, $z$ and consequently one finds the duality function for our process:
\begin{align}\label{di}
& D(x,y,z;n_{1},n_{2})= \\& \nonumber \dfrac{\bigl[ (-\frac{\sqrt{2}}{2} cos\varphi -\frac{\sqrt{2}}{2\sqrt{3}} sin\varphi)x+(\frac{\sqrt{2}}{2} cos\varphi -\frac{\sqrt{2}}{2\sqrt{3}} sin\varphi)y+\frac{\sqrt{2}}{\sqrt{3}}sin \varphi z  \bigl]^{2n_{1}}}{(2n_{1}-1)!!}\cdot \\ \nonumber &
\cdot \dfrac{\bigl[ (\frac{\sqrt{2}}{2} sin\varphi -\frac{\sqrt{2}}{2\sqrt{3}} cos\varphi)x+(-\frac{\sqrt{2}}{2} sin\varphi -\frac{\sqrt{2}}{2\sqrt{3}} cos\varphi)y+\frac{\sqrt{2}}{\sqrt{3}}cos \varphi z \bigl]^{2n_{2}}}{(2n_{2}-1)!!}.
\end{align}

\begin{oss}
Due to the arbitrary of $\varphi$, this example shows that there are infinitely many duality functions not trivially related by a multiplicative factor.
\end{oss}
The non-trivial duality relation discussed in the previous section is formalized in the following theorem.
\begin{teor}\label{original}
The process $ \mathbf{X}(t)= \bigl(  x(t), y(t), z(t) \bigl)  $  with generator $ \mathscr{L}_{3} $ in equation \eqref{elle3} is dual to the process $\mathbf{N}(t)= \bigl(  n_{1}(t), n_{2}(t) \bigl)$, with generator
\begin{align*}
L_{3}f(n_{1},n_{2}) & =n_{1}\left( n_{2}+\frac{1}{2}\right) \bigr[  f(n_{1}-1,n_{2}+1)-f(n_{1},n_{2})\bigr]  \\
& +n_{2}\left( n_{1}+\frac{1}{2}\right)  \bigr[ f(n_{1}+1,n_{2}-1)-f(n_{1},n_{2})\bigr]
\end{align*}
on duality function $D=D(x,y,z;n_{1},n_{2})$ in \eqref{di}.
\end{teor}
\dim \ \ The result is proved by an explicit computation that shows that the definition of duality in \eqref{dua} is satisfied.
\cvd

\subsection{Duality theory under change of coordinates}
\begin{figure}[hbtp]
\centering
\includegraphics[scale=0.3]{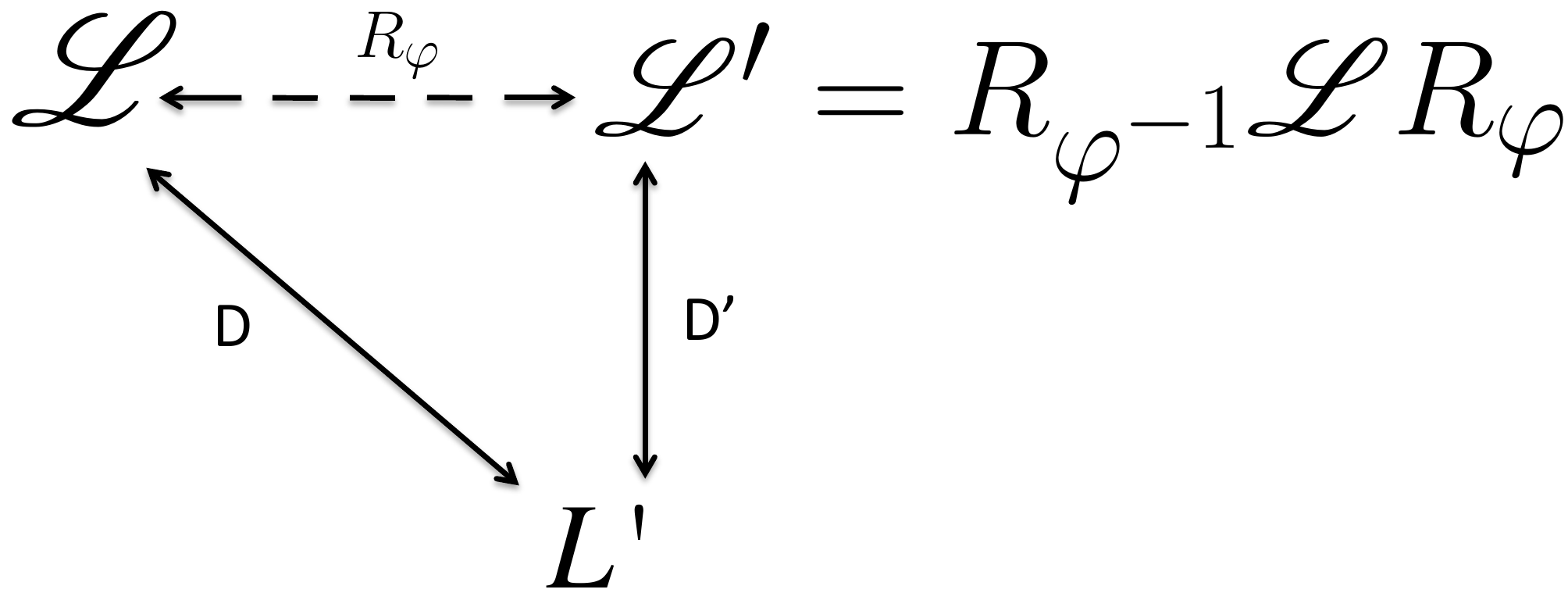} 
\caption{\small{A schematic representation of theorem \ref{change}, full arrows indicate a duality relation, while dashed map $R_{\varphi}$ is a relation between two function spaces.}}
\label{fig:change}
\end{figure}

The purpose of this last section is to extend for a generic situation the duality result found for the process that leaves energy and momentum invariant.
Starting from a known duality and after a change of coordinates, we ask about the duality properties in the new system of coordinates.
The main idea is explained in \figurename~\ref{fig:change}, where full arrows allude to a duality relation: $D'$ represents the known duality function, while $D$ is the new one. The relation given by dashed arrow $R_{\varphi}$, as formalized in theorem \ref{change}, is between the same function spaces, for which we use two different coordinate systems.

\begin{teor} [Duality and change of coordinates] \label{change}
Consider the function spaces $\mathscr{B}(\Omega)= \lbrace f:\Omega \rightarrow \mathbb{R} \rbrace $ and $\mathscr{B}(\Omega')= \lbrace f:\Omega' \rightarrow \mathbb{R} \rbrace $, where $\Omega$ and $ \Omega' $ are two generic spaces. Let 
\begin{align*}
R_{\varphi}: \mathscr{B} (\Omega') & \rightarrow \mathscr{B}(\Omega) \\ 
f  & \rightarrow \bigl( R_{\varphi}f \bigl) (x) = \bigl( f\circ \varphi \bigl)(x)
\end{align*}
the composition with function 
\begin{align*}
\varphi: \Omega & \rightarrow \Omega' \\
x & \rightarrow \varphi(x)=x',
\end{align*}
where we assume that $\varphi$ is an invertible function.\\
Consider also the two operators $\mathscr{L'}: \mathscr{B} (\Omega') \rightarrow \mathscr{B} (\Omega') $ and $L': \mathscr{B} (\Omega^{dual}) \rightarrow \mathscr{B} (\Omega^{dual})$ dual to each other with duality function $D'(x',n'): \Omega' \times \Omega^{dual} \rightarrow \mathbb{R}$, i.e.,
\begin{equation}\label{ddd}
\bigl[ \mathscr{L'}D'(\cdot,n') \bigl](x')=\bigl[ L'D'(x', \cdot) \bigl](n').
\end{equation}
Let $\mathscr{L}: \mathscr{B} (\Omega) \rightarrow \mathscr{B} (\Omega)$ be an operator related to $\mathscr{L'}$ as follows 
\begin{equation}\label{dddd}
\mathscr{L'}=R_{\varphi^{-1}}\mathscr{L}R_{\varphi}.
\end{equation}
Then, operator $\mathscr{L}$ is dual to $L'$ through $D: \Omega \times \Omega^{dual} \rightarrow \mathbb{R}$, i. e. 
\begin{equation}
\bigl[ \mathscr{L}D(\cdot,n') \bigl] (x) = \bigl[ L'D(x, \cdot) \bigl] (n')
\end{equation} 
where the duality function $D$ is
\begin{equation}\label{dprimo}
D(x,n'):=\bigl( R_{\varphi}D'(\cdot, n') \bigl) (x) = D'(\varphi(x),n').
\end{equation}
\end{teor}
\dim \ \ Consider the following chain of equalities given by \eqref{ddd}, \eqref{dddd} and \eqref{dprimo} 

\begin{equation*}
L'D(x,\cdot)=\mathscr{L'}D'(\cdot,n')=R_{\varphi^{-1}}\mathscr{L}R_{\varphi}D'(\cdot,n')
=R_{\varphi^{-1}}\mathscr{L}D(x,\cdot)
\end{equation*}
Multiplying each side for $R_{\varphi}$ we have
\begin{equation*}
\mathscr{L}D(x,\cdot)=R_{\varphi}L'D'(\cdot,n')=L'R_{\varphi}D'(\cdot,n')=L'D(x,\cdot)
\end{equation*}
\cvd

\end{document}